\documentstyle{article}
\begin{document}
\title{{Singular  Solutions of Hessian Fully Nonlinear Elliptic Equations} }
 \author{{Nikolai Nadirashvili\thanks{LATP, CMI, 39, rue F. Joliot-Curie, 13453
Marseille  FRANCE, nicolas@cmi.univ-mrs.fr},\hskip .4 cm Serge
Vl\u adu\c t\thanks{IML, Luminy, case 907, 13288 Marseille Cedex
FRANCE, vladut@iml.univ-mrs.fr} }}

\date{}
\maketitle

\def\n{\hfill\break} \def\al{\alpha} \def\be{\beta} \def\ga{\gamma} \def\Ga{\Gamma}
\def\om{\omega} \def\Om{\Omega} \def\ka{\kappa} \def\lm{\lambda} \def\Lm{\Lambda}
\def\dl{\delta} \def\Dl{\Delta} \def\vph{\varphi} \def\vep{\varepsilon} \def\th{\theta}
\def\Th{\Theta} \def\vth{\vartheta} \def\sg{\sigma} \def\Sg{\Sigma}
\def\bendproof{$\hfill \blacksquare$} \def\wendproof{$\hfill \square$}
\def\holim{\mathop{\rm holim}} \def\span{{\rm span}} \def\mod{{\rm mod}}
\def\rank{{\rm rank}} \def\bsl{{\backslash}}
\def\il{\int\limits} \def\pt{{\partial}} \def\lra{{\longrightarrow}}

\section{Introduction}
\bigskip

In this paper we study a class of fully nonlinear second-order
elliptic equations of the form
$$F(D^2u)=0\leqno(1)$$
defined in a domain of ${\bf R}^n$. Here $D^2u$ denotes the
Hessian of the function $u$. We assume that
 $F$ is a Lipschitz  function defined on an open set
$D\subset S^2({\bf R}^n)$ of the space of ${n\times n}$ symmetric
matrices  satisfying the uniform ellipticity condition, i.e. there
exists a constant $C=C(F)\ge 1$ (called an {\it ellipticity
constant\/}) such that $$C^{-1}||N||\le F(M+N)-F(M) \le C||N||\;
\leqno(2)$$ for any non-negative definite symmetric matrix $N$; if
$F\in C^1(D)$ then this condition is equivalent to
$${1\over C'}|\xi|^2\le
F_{u_{ij}}\xi_i\xi_j\le C' |\xi |^2\;,\forall\xi\in {\bf
R}^n\;.\leqno(2')$$
 Here, $u_{ij}$ denotes the partial derivative
$\pt^2 u/\pt x_i\pt x_j$. A function $u$ is called a {\it
classical\/} solution of (1) if $u\in C^2(\Om)$ and $u$ satisfies
(1).  Actually, any classical solution of (1) is a smooth
($C^{\alpha +3}$) solution, provided that $F$ is a smooth
$(C^\alpha )$ function of its arguments.

For a matrix $S \in   S^2({\bf R}^n)$  we denote by $\lambda(S)=\{
\lambda_i : \lambda_1\leq...\leq\lambda_n\}
 \in {\bf R}^n$  the (ordered) set  of eigenvalues of the matrix $S$. Equation
(1) is
 called a Hessian equation \linebreak ([T1],[T2] cf. [CNS]) if the function
 $F(S)$ depends only
  on the eigenvalues $\lambda(S)$ of the matrix $S$, i.e., if
 $$F(S)=f(\lambda(S)),$$
 for some function $f$  on ${\bf R}^n$ invariant under  permutations of
 the coordinates.

 In other words the equation (1) is called Hessian if it is invariant under
 the action of the group
 $O(n)$ on $S^2({\bf R}^n)$:
 $$\forall O\in O(n),\; F({^t O}\cdot S\cdot O)=F(S) \;.\leqno(3) $$

  If we assume that the function $F(S)$ is defined for any symmetric matrix $S$,
 i.e., $D=S^2({\bf R}^n)$ the Hessian invariance relation (3) implies the
following:

\medskip
 (a) $F$ is a smooth (real-analytic) function of its arguments if and only if $f$ is
a smooth (real-analytic) function.

\medskip
 (b) Inequalities (2) are equivalent to the inequalities
 $${\mu\over C_0} \leq { f ( \lambda_i+\mu)-f ( \lambda_i) } \leq C_0 \mu,
 \; \forall  \mu\ge 0,$$
 $\forall  i=1,...,n$, for some positive constant $C_0$.

\medskip
 (c) $F$ is a concave function if and only if $f$ is concave [Ba, CNS].

\medskip
 Well known examples of the Hessian equations are Laplace, Monge-Amp\`ere, and Special Lagrangian equations.

  \medskip We are interested also in   Isaacs equations which are uniformly elliptic but in
 general not Hessian.
 Bellman and Isaacs equations appear in the theory of controlled diffusion processes. The
 both are fully nonlinear uniformly elliptic equations of the form (1). The Bellman equation
 is concave in $D^2u \in  S^2({\bf R}^n)$ variables. However, Isaacs operators are, in  general,
 neither concave nor convex. In a simple homogeneous form the Isaacs equation can be written
 as follows:
$$F(D^2u)=\sup_b \inf_a L_{ab}u =0, \leqno (4) $$
 where $L_{ab}$ is a family of linear uniformly elliptic operators  with an ellipticity
 constant $C>0$ which depends on two parameters $a,b$.
Consider the Dirichlet problem
$$\cases{F(D^2u)=0 &in $\Om$\cr
u=\vph &on $\pt\Om\;,$\cr}\leqno(5)$$ where  $\Omega \subset {\bf
R}^n$ is a bounded domain with smooth boundary $\partial \Omega$
and $\vph$ is a continuous function on $\pt\Om$.

We are interested in the problem of existence and regularity of
solutions to  Dirichlet problem (5) for Hessian and Isaacs equations.
Dirichlet problem (5) has always a unique viscosity (weak)
solution for  fully nonlinear elliptic equations (not necessarily
Hessian equations). The viscosity solutions  satisfy the equation
(1) in a weak sense, and the best known interior regularity
([C,CC], cf. [T3]) for them is $C^{1+\epsilon }$ for some $\epsilon
> 0$. For more details see [CC,CIL]. Until recently it remained
unclear whether non-smooth viscosity solutions exist. In [NV1] we
proved the existence of viscosity solutions to the fully nonlinear
elliptic equations which are not classical in dimension 12.
Moreover, we proved in [NV2], that in 24-dimensional space the
optimal interior regularity of viscosity solutions of fully
nonlinear elliptic equations is no more than $C^{2-\delta }$. Both
papers [NV1,NV2]  use   the function
$$w={Re (q_1q_2q_3)\over |x|}, $$
where $q_i\in {\bf H},\ i=1,2,3,$ are Hamiltonian quaternions,
$x\in {\bf H}^3={\bf R}^{12}$  which is a viscosity solution in
${\bf R}^{12}$ of a uniformly elliptic equation (1) with a smooth
$F$. The proofs use some remarkable algebraic identities verified by
(the spectrum of the Hessian of)
the function $w.$ One notes also that the example by
Harvey-Lawson-Osserman [LO,HL]
 of a Lipshitz non-analytic solution to the associator (minimal surface)  equation 
strongly resembles our function. Moreover a suitable version of
 an octonion analogue [NV3] of $w$ is  reminds  the  associative
 calibration and its modifications remind coassociative and Caley
 calibrations [HL]. In our opinion these connections deserve a further study.

\medskip
The main goal of this paper is to show that the same function $w$
is a solution to a Hessian equation. Moreover the following
theorem holds

\bigskip
{\bf Theorem 1.1.} {\it For any $\delta , \; 0\leq \delta < 1$ the
function $$w/ |x|^{\delta } $$ is a viscosity solution to a
uniformly elliptic Hessian equation $(1)$ in a unit ball $B\subset
{\bf R}^{12}$.}

\bigskip
Theorem 1.1 shows that the second derivatives of  viscosity
solutions of Hessian equations (1) can blow up in an interior
point of the domain and that the optimal interior regularity of
the viscosity solutions of Hessian equations is no more than
$C^{1+\varepsilon}$, thus showing the {\em optimality} of the
result by Caffarelli-Trudinger [C,CC,T3] on the interior
$C^{1,\alpha}$-regularity of viscosity solutions of fully
nonlinear equations. Our construction provides a Lipschitz
functional $F$ in Theorem 1.1. Using a  more complicated argument
one can make $F$ smooth; we will return to this question
elsewhere. However, if we drop the invariance condition (3) we get

\bigskip
{\bf Corollary 1.1.} {\it For any $\delta , \; 0\leq \delta < 1$
the function $$w/ |x|^{\delta } $$ is a viscosity solution to a
uniformly elliptic  (not necessarily  Hessian) equation $(1)$  in
a unit ball $B\subset {\bf R}^{12}$ where $F$ is a $(C^{\infty})$
smooth functional.}

 \medskip
We show that the same function  is a viscosity solution  to
a uniformly elliptic Isaacs equation:
\bigskip

{\bf Theorem 1.2.}

{\it For any $\delta , \; 0\leq \delta < 1$   the function
$$w/ |x|^{\delta } $$ is a viscosity solution  to  a uniformly elliptic  Isaacs equation
  $(1.4)$ in a unit ball $B\subset {\bf
R}^{12}$.}

\bigskip
The question on the minimal dimension $n$ for which there exist
nontrivial homogeneous order 2 solutions of (1) remains open. We
notice that from the result of Alexsandrov [A] it follows that any
homogeneous order 2 solution of the equation (1) in ${\bf R}^3$
with a real analytic $F$ should be a quadratic polynomial. For a
smooth and less regular $F$  similar results in the dimension 3
can  be found in [HNY].

\medskip
However,we are able  reduce this dimension
 by one to 11.  Moreover the following theorem holds

\bigskip
{\bf Theorem 1.3.} {\it For any hyperplane $H\subset {\bf R}^{12}$
the function $w $ restricted to $H= {\bf R}^{11}$ is a viscosity
solution to a uniformly elliptic Hessian equation $(1)$ in a unit
ball $B\subset {\bf R}^{11}$ where $F$ is a Lipschitz functional.}

\bigskip
If we drop the invariance condition (3) we get

\bigskip
{\bf Corollary 1.2.} {\it For any hyperplane $H\subset {\bf
R}^{12}$ the function $w $  restricted to $H= {\bf R}^{11}$ is a
viscosity solution to a uniformly elliptic (not necessarily
Hessian)equation $(1)$ in a unit ball $B\subset {\bf R}^{11}$
where $F$ is a $(C^{\infty})$ smooth functional.}

\medskip\bigskip

Note, however that our technique here is {\it not} sufficient to
get singular (i.e. with unbounded second derivatives) solution in
eleven dimensions, see Remark 6.2 below.

\bigskip
 Ball $B$ in Theorem 1.1 can not be substituted
by the whole space ${\bf R}^{12}$. In fact, for any $0<\alpha <2$
there are no homogeneous order $\alpha $ solutions to the fully
nonlinear elliptic equation (1) defined in ${\bf R}^n\setminus \{
0 \}  $, [NY]; the essence of the difference with the local
problem is that in the case of homogeneous solution defined in
${\bf R}^n\setminus \{ 0 \}$ one deals simultaneously with two
singularities of the solution: one at the origin and another at
the infinity. In the local problem the structure of singularities
of solutions is quite different, even in dimension 2, the function
$u=|x|^{\alpha } ,\; 0<\alpha < 1, \;   x\in B^o$, where $B^o$ is
a punctured ball in ${\bf R}^n, \;  n\geq 2, \; B^o= \{ x\in {\bf
R}^n, 0< |x|<1 \} $, is a solution to the uniformly elliptic
Hessian equation in $ B^o $ (notice that $u$  {\it is not } a
viscosity solution of any elliptic equation on the whole ball
$B$).

  We study also the possible singularity of solutions of Hessian equations
  defined in a neighborhood of a point. We prove the following general
result:

\bigskip
 {\bf Theorem 1.4.} {\it Let $u$ be a viscosity solution of a uniformly elliptic Hessian
  equation in a punctured ball $B^o\subset {\bf R}^n$. Assume that $u\in C^0(B)$.
 Then $u=v+l+o(|x|^{1+\varepsilon }) $, where
 $v$ is a monotone function of the radius, $v(x)=v(|x|)$, $v\in C^{\varepsilon}(B)$,
  where $\epsilon >0$ depends on the ellipticity constant of the equation,
 and $l$ is a linear
function.}

\bigskip
As an immediate consequence of the theorem we have

\bigskip
 {\bf Corollary 1.3.} { \it Let $u$ be a homogeneous order $\alpha ,\;  0<\alpha  < 1$
 solution of a uniformly elliptic Hessian equation in
 a punctured ball $B^o\subset {\bf R}^n$. Then $u=c|x|^{\alpha }$.}

\medskip
 The rest of the paper is organized as follows: in Section
2 we give a sufficient condition for validity of Theorem 1.1, we
verify it in Section 3 for $\delta=0$   and then in Section 4 for
any $1>\delta\ge 0$. Section 5 is devoted to a proof of Theorem
1.2, Section 6 proves Theorem 1.3, and Section 7 contains a proof of Theorem 1.4.

\medskip
{\em Acknowledgement.} The authors would like to thank L.
Caffarelli who posed the question leading to the present work.

\medskip
Since the proof of Theorem 1.1 in Sections 3 and 4 is somewhat involved and utilize 
computer (MAPLE) computations, we give here an account of its logical structure and 
its principal  points. First of all, the criterion of ellipticity in Section 2 reduces Theorem 1.1 
for $\delta =0$ to the uniform hyperbolicity of  $Hess(P)(a)-{^t O }\cdot Hess(P)(b)\cdot O$
for a pair $a\neq b$ of unit vectors and an orthogonal matrix $O$. A classical result by H. Weyl
 on the eigenvalues of the diference of two symmetric matrices reduces this to the uniform 
hyperbolicity of the difference $\lambda(Hess(P)(a))-\lambda( Hess(P)(b))$.
 Recall then  [NV1, Section 3] that  the characteristic polynomial $CH(P,a)(T)$
of the Hessian $Hess(P)(a)$ of the cubic form $P$ has  for $a\in S^{11}_1$ the following form:
$$   CH(P,a)(T)=(T^3-T+2m(a))(T^3-T-2m(a))(T^3-T+2P(a))^2, \; $$ where  $m(a)\ge |P(a)|$
which permits to conclude that the structure of the (ordred) spectrum is  as follows
$$\mu_1=\mu'_1\ge \lambda_1\ge \lambda_2\ge \lambda_3\ge \mu_2=\mu'_2\ge -\lambda_3\ge-
 \lambda_2\ge - \lambda_1\ge \lambda_3\ge\mu_3=\mu'_3$$
where $\mu_1\ge\mu_2\ge\mu_3$ are the roots of $(T^3-T+2P(a))$,  and 
$\lambda_1\ge \lambda_2\ge \lambda_3\ge \  -\lambda_3\ge-
 \lambda_2\ge - \lambda_1$ are those of $(T^3-T+2m(a))(T^3-T-2m(a))$. 
The argument of Section 3 is based on the calculation of the (shifted)
 characteristic polynomial $CH(w,a)(T-P(a))$ of the full Hessian $Hess(w)(a)$ which is possible
 thanks to an action of the group $Sp(1)\times Sp(1)\times Sp(1)$ which does not change 
this polynomial. This action permits to bring the matrix  $Hess(w)(a)$ to a simple block form and 
gives using a MAPLE caluculation an explicit formula for $CH(w,a)(T-P(a))$:
$$CH(w,a)(T-P(a))=P_6(a,T)(T^3-T+2P(a))^2$$
for a certain explicit polynomial $P_6(a,T);$  in fact $P_6(a,T)$ is  the (shifted)
 characteristic polynomial of $Hess(w_6)(a')$ for a 6-dimensional version of $w$ and an appropriate 
6-dimensional unit vector $a'$.  The crucial point then is that the spectrum
 in this case is not so different from that of $Hess(P)(a)$. 	In fact, one has for this ordered spectrum:
$$\mu_1=\mu'_1\ge \lambda'_1\ge \lambda'_2\ge \lambda'_3\ge \mu_2=\mu'_2\ge \lambda'_4\ge
 \lambda'_5\ge  \lambda'_6\ge \lambda_3\ge\mu_3=\mu'_3$$
where $\lambda'_1\ge \lambda'_2\ge \lambda'_3\ge  \lambda'_4\ge
 \lambda'_5\ge  \lambda'_6$ are the roots of $P_6(a,T)$. To prove this inequalities one verifies 
it for specific points $a$ and then explicitly calculates (using MAPLE)
 the resultant which (miraculously)  vanishes nowhere and thus gives the necessary inequalities.
This garanties the exact formula for the equal  6th and  7th eigenvalues which permits to get 
the necessary uniform hyperbolicity of the difference $\lambda(Hess(P)(a))-\lambda( Hess(P)(b))$.

In Section 4 we generalize this agument to any $\delta\in ]0,1[$. In this situation we need 
the uniform hyperbolicity of  $Hess(P)(a)-K{^t O }\cdot Hess(P)(b)\cdot O$
for a pair $a\neq b$ of unit vectors, any orthogonal matrix $O$ and any positive constant $K$, 
which follows from that of $\lambda(Hess(P)(a))-K\lambda( Hess(P)(b))$. We begin with
 the uniform hyperbolicity of the  difference $(\mu_1(a),\mu_2(a),\mu_3(a))-K(\mu_1(b),\mu_2(b),\mu_3(b))$ 
 which is rather elementary since there are simple trigonometric formulas for $\mu_i$. Unfortunately,
 the position of $\mu_2$ in the ordered spectrum of  $Hess(P)(a)$ is not   fixed anymore, which
follows from an explicit calculation  of $CH(w_{\delta},a)(T-(1+\delta)P(a))$  together with some 
resutant calculations similar (but more involved) to those  in Section 3. However, the position
of the double value $\mu_2=\mu'_2$ varies from (5,6) to (7,8) and an argument using the oddness 
of $w_{\delta}$ permits to deduce the uniform hyperbolicity of 
$\lambda(Hess(P)(a))-K\lambda( Hess(P)(b))$ from that of
 $(\mu_1(a),\mu_2(a),\mu_3(a))-K(\mu_1(b),\mu_2(b),\mu_3(b))$ which finishes the proof 
of Theorem 1.1.

\section{Ellipticity }

\bigskip

   Let $w$ be a homogeneous function
of order $2-\delta ,\  0\leq \delta <1 $, defined on a unit ball
$B =B_1\subset {\bf R}^n$ and smooth in $B \setminus\{0\}$. Then
the Hessian of $w$ is homogeneous of order $(-\delta)$. Define the
map
$$\Lambda :B  \longrightarrow \lambda (D^2w) \in {\bf R}^n\; .$$

Let $K\subset {\bf R}^n$ be an open convex cone, such that
$$ \{ x\in {\bf R}^n : x_i\geq 0,\ i=1,...,
n\}  \subset K .$$

Set
$$L:={\bf R}^n\setminus (K\cup -K).$$

  We say that a set  $E\subset {\bf R}^n$ satisfy $K$-cone condition if
  $\;  (a-b)\in L$
for any $ a,b\in E.$

Let $S_n$ be the group  of permutations of   $\{ 1,...,n\}$. For
any $\sigma \in S_n$, we denote by $T_{\sigma}$ the linear
transformation of ${\bf R}^n$ given by
 $x_i \mapsto x_{\sigma(i)}, \; i=1,...,n.$

\bigskip\noindent {\bf  Lemma 2.1.} {\it  Assume that
$$M:=\bigcup_{\sigma \in S_{n} }\   T_{\sigma }\Lambda (B)\subset {\bf R}^n $$
  satisfies the $K$-cone condition. If $\delta> 0$
we assume additionally that $w$ changes sign in $B$. Then $w$ is a
viscosity solution in $B$ of a uniformly elliptic Hessian equation
$(1)$.}
\bigskip\noindent

{\em   Proof  }. Let us choose in the space ${\bf R}^n$ an
orthogonal coordinate system $z_1,\dots,z_{n-1},s,$  such that
$s=x_1+...+x_n$ . Let $\pi : {\bf R}^n\to Z$ be the orthogonal
projection of ${\bf R}^n$ onto the $z$-space. Let $K^\ast$ denote
the adjoint cone of $K$, that is, $K^\ast = \{b\in {\bf R}^n:
b\cdot c \ge 0 \ for\ all\ c \in K\} $. Notice that $a \in L $
implies $a\cdot b =0$ for some $b \in K^\ast$.  We represent the
boundary of the cone $K$ as the graph of a Lipschitz function
$s=e(z)$, with $e(0)=0$, function $e$ is smooth outside the
origin:
$$e(z)\ =\inf\{ c:\ (z+cs)\in K \} .$$

Set $m=\pi \bigl (M )$. We prove that $M$ is a graph of a
Lipschitz  function on $m$,
$$M =\{z\in m:s=g(z)\}\;.$$
Let $a,\hat a \in M ,  a =(z,s),\hat a = (\hat z,\hat s)$. Since
$a-\hat a \in L$, we have $$-e(z-\hat z) \le \hat s - s \le e(z-
\hat z).$$ Since $e(0)=0, g(z):=s$ is single-valued. Also
$$|g(z)-g(\hat z)|=|s-\hat s|\le |e(z-\hat z)|\le C|z-\hat z|.$$

The function $g$ has an extension $\widetilde{g}$ from the set
$m$ to ${\bf R}^{n-1}$ such that $\widetilde{g}$ is a
Lipschitz function and the graph of $\widetilde{g}$ satisfies the
$K$-cone condition. One can define such extension $\widetilde{g}$
simply by the formula
$$\widetilde{g}(z):=\inf_{w\in m } \bigl\{ g(w)\ +  e(z-w)
\bigr\}\;.$$

To show that this formula works let $(z, \tilde g(z)),  (\hat
z,\tilde g(\hat z))$ lie in the graph $\tilde g$. We must show
$$-e(z-\hat z)\ \le \ \tilde g(z)-\tilde g(\hat z)\ \le e(z-\hat z).$$
Now
$$\tilde g(\hat z)\ =\ g(w)+e(\hat z-w)$$
for some $w\in m$. Thus
$$\tilde g(z)-\tilde g(\hat z)\ \le \ g(w)+e(z-w)-(g(w)+e(\hat z-w))\ \le \
e(z-\hat z),$$ since $e(a+b)\ \le \ e(a)+e(b)$, as $e(\cdot )$ is
convex, homogenous. Similarly
$$\tilde g(z)-\tilde g(\hat z)  \ge  -e(z-\hat z).$$
\medskip
Let us set
$$
f':= s - \widetilde{g}(z).
$$
Since the level surface of the function  $ f'$ satisfies $K$-cone
condition it follows that $\nabla f \in K^*$ a. e.  where  $K^*$ is the
adjoint cone to $K$.
Moreover   the function $w$ satisfies the equation
$$
 f'(\lambda (D^2 w))= 0.
$$
on $B \setminus \{0\}$.

Set
$$f=\sum_{\sigma \in {S_n}} f'(\sigma (x)).$$
Then $f$ is a Lipschitz  function invariant under the action of
the group $S_n$ and satisfies the equation
$$
 f(\lambda (D^2 w))= 0.
$$
on $B \setminus \{0\}$.

\medskip
We show now that $w$ is a viscosity solution of (1) on the whole
ball $B$.

Assume first that $\delta =0$. Let $p(x),\ x\in B$ be a quadratic
form such that $p\leq w$ on $B$. We choose any quadratic form
$p'(x)$ such that $p\leq p'\leq w$ and there is a point $x'\neq 0$
at which $p'(x')= w(x')$. Then it follows that $F(p)\leq F(p')\leq
0$. Consequently for any quadratic form $p(x)$ from the inequality
$p\leq w$ ($p\geq w$) it follows that $F(p)\leq 0$ ($F(p)\geq 0$).
This implies that $w$ is a viscosity solution of (1) in $B$ (see
Proposition 2.4 in [CC]).

If $0< \delta < 1$ then for any smooth function $p$ in  $B$ the
function $w-p$ changes sign in any neighborhood of $0$. Hence, by
the same proposition in [CC], it follows that $w$ is a viscosity
solution of (1) in $B$.
\bigskip

\section{Non-classical solution}
This section is devoted to a proof of Theorem 1.1 in the case of
$\delta=0$ i.e. for a non-classical, but not singular, solution.

We define the cubic form $P$ which is used to construct our
non-classical and singular solutions. Let $ X=(r,s,t)\in {\bf
R}^{12}$ be a variable vector with $r,s,$ and $t\in {\bf R}^4.$
For any $ t=(t_0,t_1,t_2,t_3)\in {\bf R}^4$ we denote  by $
qt=t_0+t_1\cdot i+t_2\cdot j+t_3\cdot k\in {\bf H} $ (Hamilton
quaternions).

Define the cubic form $P=P(X)=P(r,s,t) $ as follows
$$P(r,s,t)=Re(qr\cdot qs\cdot qt)=r_0s_0t_0-r_0s_1t_1-r_0s_2t_2-r_0s_3t_3$$
$$-r_1s_0t_1- r_1s_1t_0-r_1s_2t_3+r_1s_3t_2-r_2s_0t_2+r_2s_1t_3-r_2s_2t_0-r_2s_3t_1 $$
$$ -r_3s_0t_3- r_3s_1t_2+r_3s_2t_1-r_3s_3t_0;$$
and denote
$$ w(X)=P(X)/ |X|  . $$
Note that by definition one has $  | P(X)|\le  {|X|^3\over 3 \sqrt
3} ,$ since $$ |P(r,s,t)|\le |r|\cdot| s|\cdot |t|\le {\left({
r^2+ s^2+ t^2\over 3}\right)}^{3/2}.$$ In particular for $X\in
S_1^{11} $ one has
 $ |P(X)|= |w(X)|\le {1\over 3 \sqrt
3} .\;$
 For $a\in {\bf R}^{12}-\{0\} $ we denote by $H(a)$ the Hessian $D^2 w(a).$

\bigskip
 {\bf Proposition 3.1.} {\em Let  $a\neq b\in S_1^{11} $  and
 let $O\in {\hbox {O}}({12} )$ be an orthogonal matrix s.t.
 $H(a,b,O):=H(a)- {^tO}\cdot H(b)\cdot O\neq 0$.
 Denote $ \Lambda_1\ge\Lambda_2\ge
 \ldots\ge\Lambda_{12}$  the eigenvalues of the matrix
 $H(a,b,O).$
  Then}

$${1\over 26} \le {\Lambda_1\over -\Lambda_{12}}\le 26.$$
 \bigskip

 We need the following property of the eigenvalues $ \lambda_1\ge\lambda_2\ge
 \ldots\ge\lambda_{n}$ of real symmetric matrices  of order $n$:

\medskip
{\bf Property 3.1.} {\em Let $ A,B$  be two real symmetric
matrices with the eigenvalues $
\lambda_1\ge\lambda_2\ge\ldots\ge\lambda_{n} $ and $
\lambda'_1\ge\lambda'_2\ge\ldots\ge\lambda'_{n} $ respectively.
Then for the eigenvalues $
\Lambda_1\ge\Lambda_2\ge\ldots\ge\Lambda_{n} $ of the matrix $A+B$
we have
$$ \Lambda_1\ge\lambda_i+\lambda'_j, \;\;\Lambda_n\le\lambda_i+\lambda'_j$$
whenever} $i+j=n .$

\medskip
This is a classical result by Hermann  Weyl [We], cf. [Fu], p.
211.

 We will use this result in the form  which follows  (replace $B$
 by $-B$ in Property 3.1):

  \bigskip
  {\em Let $ A,B$  be two real symmetric  matrices
with the eigenvalues $
\lambda_1\ge\lambda_2\ge\ldots\ge\lambda_{n} $ and $
\lambda'_1\ge\lambda'_2\ge\ldots\ge\lambda'_{n} $ respectively.
Then for the eigenvalues $
\Lambda_1\ge\Lambda_2\ge\ldots\ge\Lambda_{n} $ of the matrix $A-B$
we have}
$$ \Lambda_1\ge\max_{i=1,\cdots, n}(\lambda_i-\lambda'_i),
 \;\;\Lambda_n\le\min_{i=1,\cdots, n}(\lambda_i-\lambda'_i).$$

\bigskip
 {\bf Main Lemma 3.1. } {\em  Let $ A:=H(a),$
$B:={^tO}\cdot H(b)\cdot O.$

$(i)$ If $P(a)-P(b)\ge 0 $ then }${\hbox
{Tr}}(B-A)=15(P(a)-P(b))\le 15\Lambda_1; $

$(ii)$ {\em If $P(a)-P(b)\le 0 $ then} ${\hbox
{Tr}}(B-A)=15(P(a)-P(b))\ge 15\Lambda_{12}. $

\bigskip
 {\em Proof of Proposition 3.1.} We consider only the case
 ${\hbox {Tr}}(A-B)=15(P(b)-P(a))\ge 0$, the proof in the other
 case being symmetric. Since ${\hbox {Tr}}(A-B)=
 \Lambda_1+\Lambda_2+\ldots+\Lambda_{12}\ge 0$ one gets
 $11\Lambda_1\ge -\Lambda_{12}.$ On the other hand,
 $$-15\Lambda_{12}\ge {\hbox {Tr}}(A-B)=
\Lambda_1+\Lambda_2 +\ldots+\Lambda_{12}$$ implies
$$-26\Lambda_{12}\ge
-15\Lambda_1-\Lambda_2-\Lambda_3-\ldots-\Lambda_{12}\ge
\Lambda_{1}$$
 which finishes the proof.

\medskip
To prove Main Lemma we need two  lemmas which constitute our
principal technical tool. We postpone their proof until the end of
the section.
 \medskip

{\bf Lemma 3.2. } {\em  Let $a=(r,s,t)\in S_1^{11};$   define
$$
W=W(a)=P(a),\; m=m(a) = |s|  ,\; n=n(a)= |t|  .$$

 Then the characteristic polynomial of the
matrix $ A:=H(a)$ is given by

$$ P_A(T)=P_1(T)^2\cdot P_2(T)$$
where
$$P_1(T)=T^3+3WT^2+3W^2T-T+W+W^3,$$
$$P_2(T)=T^6+9WT^5+(21W^2+3L-2)T^4+2W(7W^2+3L-4)T^3+
$$
$$(1-6W^2-9W^4
-3L+9M)T^2-(15W^4+6W^2L-4W^2-6L+1)WT$$
$$-
5W^6-3LW^4+4W^4-3(3M+L)W^2+W^2-M$$ with}
$L:=L(m,n)=m^2+n^2-n^2m^2-n^4-m^4\in [M,{1\over 3}] ,$

$M:=M(m,n)=m^2n^2(1-n^2-m^2)\in [W^2,{1\over 27}] .$
  \bigskip

{\bf Lemma 3.3.} {\em  Let $a=(r,s,t)\in S_1^{11}, \;\; A=H(a) .$
Let $\mu_1\ge \mu_2\ge\mu_3 $ be the roots of $P_1(T)$, $\nu_1\ge
\nu_2\ge\ldots\ge \nu_6$ be the roots of $P_2(T)$. Then}
$$\mu_1\ge\nu_1\ge \nu_2\ge\nu_3\ge\mu_2\ge\nu_4\ge\nu_5\ge\nu_6\ge\mu_3.$$

  \medskip
{\bf Corollary 3.1.} {\em  Let $a=(r,s,t)\in S_1^{11}
 .$ Let $\lambda_1\ge
\lambda_2\ge\ldots\ge\lambda_{12} $ be the eigenvalues  of
$A=H(a)$. Then}
$$\lambda_6=\lambda_7=
{2\over \sqrt 3}\cos\left({\arccos(3\sqrt 3 P(a))+\pi\over 3
}\right)-P(a).$$

{\em Proof of Corollary.} By Lemmas 3.1 and 3.2
$\lambda_6=\lambda_7=\mu_2.$ One easily verifies that
 $Q_1(X):=P_1(X-W)=X^3 -X+2W.$ If we set $X=2\cos(\beta)/\sqrt 3,\; 3\sqrt 3
  W=\cos(\alpha)$ we get $\cos(3\beta)=\cos(\alpha) $  which implies $$\mu_1=
 {2\over \sqrt 3}\cos\left({\arccos(3\sqrt 3 W)-\pi\over 3
}\right)-W,  \mu_2={2\over \sqrt 3}
 \cos\left({\arccos(3\sqrt 3 W)+\pi\over 3 }\right)-W,$$ $$\;\mu_3=
 {2\over \sqrt 3}\cos\left({\arccos(3\sqrt 3 W)+3\pi\over 3
}\right)-W. $$

\medskip
 {\em Proof of Main Lemma 3.1.} Let
$W=P(a),$ $W'=P(b)$   and $W-W'\ge 0 .$ By Property 3.1.
$$\Lambda_1\ge\lambda_6(A)-\lambda_6(B)={2\over \sqrt 3}\left(
 \cos\left({\arccos(3\sqrt 3 W)+\pi\over 3 }\right)-
 \cos\left({\arccos(3\sqrt 3 W')+\pi\over 3 }\right)\right)-W
 +W'.$$ Since
 $ { \cos\left({\arccos(3\sqrt 3 W)+\pi\over 3
}\right) }\ge { \sqrt 3 |W|} $ and $ { \cos\left( {\arccos(3\sqrt
3 W')+\pi\over 3 }\right) }\ge { \sqrt 3|W'|}$   we get the
conclusion. The case $P(a)-P(b)\le 0 $ is symmetric.

\bigskip
 {\em Proof of  Lemma 3.2.} Note that the function $w$ is invariant
 under the action of the group ${\hbox {Sp}}_{1}
  \times{\hbox {Sp}}_{1} \times{\hbox {Sp}}_{1}  $ by conjugation
  on each factor, i.e.
  $$ (g_1,g_2,g_3): (r,s,t) \mapsto (g_1rg_1^{-1},g_2sg_2^{-1},g_3tg_3^{-1}) $$
  for $g_1,g_2,g_3\in {\hbox {Sp}}_{1}=\{q\in {\bf H} : \;
  |q|=1\},$ and hence the spectrum $Sp(H(a)) $ is
  invariant under this action as well.

  Applying this action one can suppose that $r_2=r_3=s_2=s_3=t_2=t_3=0,
  $ i.e. that $(r,s,t)\in {\bf C}^3\subset {\bf H}^3.$ In this case
the matrix $A=H(a)$ becomes a block matrix
  $$A=
  \left(%
\begin{array}{cc}
  A_6&0  \\
   0&M_6  \\

\end{array}%
\right)$$ where $A_6=D^2w_6(a')$ is the Hessian of the function 
$$
w_6(a')={P_6(a')\over |a'|}={Re(cr\cdot cs\cdot ct)\over |a'|}={r_0s_0t_0-r_0s_1t_1-r_1s_0t_1-r_1s_1t_0\over\sqrt{r_0^2+s_0^2+t_0^2+r_1^2+s_1^2+t_1^2}},
$$ 
$a'=(cr, cs,ct)=(r_0+r_1 i,s_0+s_1 i,t_0+t_1 i)\in {\bf C}^3,$  and $M_6$ is
the following matrix:
$$M_6=
  \left(%
\begin{array}{cccccc}
  -W&0&-t_0&-t_1& -s_0& s_1 \\
   0& -W& t_1& -t_0& -s_1&-s_0 \\
  -t_0& t_1& -W& 0& -r_0& -r_1 \\
   -t_1& -t_0& 0& -W& r_1& -r_0\\
   -s_0&-s_1& -r_0& r_1& -W& 0\\
   s_1& -s_0& -r_1& -r_0& 0& -W\\
\end{array}%
\right).$$  A direct calculation shows that the characteristic
polynomial of   $$N_6=M_6+W\cdot I_{6}=
  \left(%
\begin{array}{cccccc}
  0&0&-t_0&-t_1& -s_0& s_1 \\
   0& 0& t_1& -t_0& -s_1&-s_0 \\
  -t_0& t_1& 0& 0& -r_0& -r_1 \\
   -t_1& -t_0& 0& 0& r_1& -r_0\\
   -s_0&-s_1& -r_0& r_1& 0& 0\\
   s_1& -s_0& -r_1& -r_0& 0& 0\\
\end{array}%
\right).$$ is given by $$ P_{N_6}(X)= (X^3 -X+2W)^2$$ (one uses
that $ |a| ^2= {|a'|} ^2= |r| ^2+ |s| ^2+ |t| ^2=1$) which gives
the formula for the first factor. To caculate the characteristic
polynomial of $A_6$ one notes an action of the group
$$T^2=S^1\times S^1=\{ (u_1, u_2,u_3)\in {\bf C}^3\;: u_1=
  u_2 = u_3 =1, u_1 u_2u_3=1\}$$ on ${\bf C}^3$
respecting $w_6 $:
$$ (u_1,u_2,u_3): (r,s,t) \mapsto (u_1r,u_2s,u_3t) .$$
This action permits to suppose that $s_1=t_1=0$,  $ s',t'\in {\bf
R}^+$ and thus $s'=s_0=m,t'=t_0=n,\; W= P(r,s,t)=r_0 m n $.
 Applying MAPLE   one gets the characteristic
 polynomial $ P_2(T)$.
 
One notes also that in
this case a direct calculation gives for $A_6=(N_{ij})$:
$$
  N_{11}=(3r_0^2-3)W,\; N_{12}=(3Wr_0-mt_0)r_1,\; N_{13}=n(1-r_0^2-m^2)+3Wr_0m
   ,\; N_{14}= r_0 n r_1,\; $$
   $$N_{15}=m(1-r_0^2-n^2)+3r_0nW,\; N_{16}= r_0 m r_1,
\; N_{21}=(3Wr_0-m n)r_1,\; N_{22}=3W(r_1^2-1),\;$$
$$ N_{23}= (3Ws_0-m n)r_1,
\; N_{24}=n(r_1^2-1),\; N_{25}=(3Wn r_0 m)r_1,\;
N_{26}=m(r_1^2-1),
$$
 $$ N_{31}= (1-r_0^2-m^2)n+3r_0mW,\; N_{32}=(3mW-r_0n)r_1,\;
 N_{33}=(3m^2-3)W,\; N_{34}=
 m nr_1, $$
 $$ N_{35}=(1-m^2-n^2)r_0+3mt_0W,\; N_{36}= (m^2-1)r_1,
  \; N_{41}= r_0 n r_1,\; N_{42}=  (r_1^2-1)n,\;N_{43} =m n r_1, $$
   $$ N_{44}=-W,\; N_{45}=(n^2-1)r_1,\; N_{46}= -r_0,
  \; N_{51}=(1-r_0^2-n^2)m+3r_0nW ,\;N_{52}= (3n W-m r_0)r_1,$$
  $$ N_{53}=
(1-m^2-n^2)r_0+3m n W,\; N_{54}=  (n^2-1)r_1,\; N_{55}=
(3n^2-3)W,\; N_{56}=  m nr_1,$$
$$N_{61}=m r_0 r_1, \;
N_{62}= (r_1^2-1)m, \; N_{63}= (m^2-1)r_1,
  \; N_{64}= -r_0, \; N_{65}= m n r_1, \; N_{66}= -W $$
which permits a human (albeit very tedious) calculation of the polynomial.

\smallskip
 Note that the caracteristic polynomial
 $Q_2(X)=P_2(X-W)$ of $A_6+W\cdot I_6 $ equals
$$Q_2(X)=X^6+3WX^5-(9W^2-3L+2)X^4
-6WLX^3+
 (6W^2 -3L +9M+1)X^2$$ $$-3(6M-4L+1)WX+
 3W^2-12LW^2-M.$$
In fact, one can directly apply the MAPLE directive $$P2:=sort(factor(simplify(charpoly(hessian(w_6,v),S))),S);$$
for the coordinate vector $v$, but in this case the calculation takes about  a minute, 100 MB of space
 (and the result need many dozens lines to be written), while the same directive
 applied to the case with two zero coordinates gives the result in less than a second.
\medskip

{\em Remark 3.1.} Since $ |a| ^2= {|a'|} ^2= |r| ^2+ |s| ^2+ |t|
^2= r_0 ^2+r_1 ^2+m ^2+n ^2 =1$ one gets  $r_0 ^2+m ^2+n ^2 \le 1$
and an application
$$ \Phi: S^{11}_1\longrightarrow  \bar B^{3}_{++},  \;\; a=(r,s,t)\mapsto
 \Phi(a):=(r_0,m,n)=({W\over mn},m,n) $$
where $\bar B^{3}_{++}=\bar B^{3}_1 \bigcap \{m\ge 0,n\ge 0\}.$

\bigskip
{\em Proof of  Lemma 3.3.} Let $\mu'_i=\mu_i+W$, $\nu'_j=\nu_i+W$
for  $i=1,2,3,\;j=1,\ldots,6;$ be the roots of $Q_1(X) $ and
$Q_2(X)$, respectively. We have to show that
$$\mu'_1\ge\nu'_1\ge \nu'_2\ge\nu'_3\ge\mu'_2\ge\nu'_4\ge\nu'_5\ge\nu'_6\ge\mu'_3.$$
One notes that $\mu'_i(W)=\mu'_i(-W)$, $\nu'_i(W)=\nu'_i(-W)$.
Therefore we can suppose w.r.g. that $W\ge 0$. For $n=0 $  we have
$ W=mnr_0=0$ and
$$Q_2(X)=
X^6  - 2 X^4  + 3 m  X^4  - 3 m  X^4  + X^2  - 3 m  X^2  + 3 m
X^2=$$ $$X^2  (X - 1) (X + 1) (X^2  - 3 m^4  - 1 + 3 m^2 )
$$ $$ X^3-X-2W=X^3-X=X  (X - 1) (X + 1).$$
Thus $\mu'_1=\nu'_1=1$, $\nu'_2=\sqrt{  1 - 3 m^2+ 3 m^4  } \in
(0, 1]$, $\nu'_3=\nu'_4=\mu'_2=0$, $\nu'_5=-\sqrt{  1 - 3 m^2+ 3
m^4 } \in [-1, 0)$, $\mu'_3=-1,$ and the inequalities take place.
Symmetrically this is true for $m=0$ as well. We can suppose thus
that $m\neq 0,\;n\neq 0. $

We suppose  then that  $ r_0^2+m^2+n^2 \neq 1;$  without loss  one
supposes also  $(m,n,r_0 )\in B_1^3\bigcap {\bf R}_+^3 .$ We begin
with a particular choice: $m=n=r_0=1/2,\; W=1/8. $ For that choice
easy brute force calculations show that $ \mu'_1\in [0.83, 0.84]$,
$ \mu'_2\in [0.26, 0.27]$,  $ \mu'_3\in [-1.11,-1.1]$,
$\nu'_1\in[0.7,0.71]$, $\nu'_2\in[0.54,0.55]$,
$\nu'_3\in[0.42,0.43]$,$\nu'_4\in[-0.39,-0.38]$,
$\nu'_5\in[-0.71,0.7]$, $\nu'_6\in[-0.96,0.95]$ and the
inequalities hold. Then we consider the resultant $R=R(m,n,r_0)$
of the polynomials $Q_2(X)$ and $X^3-X+2W; $ a brute force (MAPLE)
calculations give
$$R=16(-n^2m^2+W^2+n^2m^4+n^4m^2)^3(27W^2+4)(1-27W^2)= $$
$$
16n^2m^2(-1+r_0^2+m^2+n^2)^3(27W^2+4)(1-27W^2)< 0$$ since  the
condition  $ W^2 =1/27$  implies $ r_0^2+m^2+n^2=1.$ For any
$(m,n,r_0 )\in B_1^3\bigcap {\bf R}_+^3 $  there is a line segment
joining it to the triple $(1/2,1/2,1/2 )$,  the set $ B_1^3\bigcap
{\bf R }_+^3 $ being convex. The value of $R(m,n,r_0)$ on the
whole segment is strictly negative and thus the order of the roots
at $( m, n, r_0 )$ is the same as at $(1/2,1/2,1/2 )$ which
finishes the proof of the inequalities for $ r_0^2+m^2+n^2 \neq
1$. Let finally $m^2+n^2+r_0^2=1.$ Then easy
 brute force calculations show that
$$Q_2(X)=(X^3-X+2W)(X^3+3WX^2-9W^2X-X+3LX +W-6WL).$$
  Thus by continuity we get
$\lambda_1=\lambda_2=\lambda_3=\mu_1,
\;\lambda_{6}=\lambda_{7}=\mu_2,
\lambda_{10}=\lambda_{11}=\lambda_{12}=\mu_3 $ which is sufficient
to conclude.

\bigskip
{\em Remark 3.2.} We use extensively MAPLE calculations in
Sections 3 and 4. These calculations concern algebraic identities,
do not use any approximation and are thus completely  rigorous.
Besides, all of them need only few seconds on a modest laptop.

\bigskip
 Proposition 3.1 and Lemma 2.1 give a proof of Theorem 1.1 in the case of
  $ \delta=0\; .$ Indeed,  we set $K$ to be the dual cone  $ K:=K_\lambda^{\ast}$ where
  $$  K_\lambda=
\{( \lambda_1,...  ,\lambda_{n})\in [C/\lambda, C\lambda]: {\hbox{ for some}}\; C>0 \;\}
$$ with $ n=12, \;  \lambda=26.$ Then Proposition 3.1 gives  the $K-$cone condition in
 Lemma 2.1 on   $T_{\sigma_0 }\Lambda (B)$ for ${\sigma_0=id \in S_{12} }$
which implies the same  condition on the whole
 $M=\bigcup_{\sigma \in S_{n} }\   T_{\sigma }\Lambda (B) $ as well.

 \bigskip
{\em Remark 3.3.}  The ellipticity constant $C$ of thus obtained
functional $F$ verifies
 $C\le 4\cdot 26^2\sqrt{12}<10^5 $ (cf. [NV1, Lemma 2.2]).

 \section{Singular solutions}

In this section we prove Theorem 1.1 for any $ \delta\in [0,1)\;
.$ For this it is sufficient to show  by Lemma 2.1 that the
ellipticity condition (the $K-$cone condition) valid for the
function $w$ remains to hold for the function $w_{\delta}(X):=w(X)
|X| ^{-\delta}.$

 For $a\in {\bf R}^{12}-\{0\} $ we denote
by $H_{\delta}(a)$ the Hessian $D^2w_{\delta}(a).$ The following
result is sufficient to prove Theorem 1.1:

 \bigskip
{\bf Proposition 4.1.}
   {\em  Let $0\le {\delta}<1. $ Then  for any
   $a\neq b\in {\bf R}^{12}-\{0\}  $   and any
   orthogonal matrix $O\in O{(12)} $  with
 $H_{\delta}(a,b,O):=H_{\delta}(a)- {^tO}\cdot H_{\delta}(b)\cdot O\neq
 0$ the eigenvalues
 $ \Lambda_1\ge\Lambda_2\ge
 \ldots\ge\Lambda_{12}$  of
 $H_{\delta}(a,b,O)$
  verify}

$${ 1 \over C_{\delta}}={ 1-\delta \over 26+3\delta-\delta^2} \le {\Lambda_1\over -\Lambda_{12}}
\le { 26+3\delta-\delta^2 \over 1-\delta}=:C_{\delta}\; .$$
 \bigskip

{\em Proof.} We can suppose  without loss that $ |a| \le  |b|$,
moreover, by homogeneity we can suppose that $a\in S_1^{11}$ and
thus  $ |b| \ge 1.$ Let $\bar b:={b/|b|}\in S_1^{11}$ then
$D^2w_{\delta}(b)=D^2w_{\delta}(\bar b)|b|^{-\delta} .$ One needs
then the following  result for the points $a,\bar b \in S_1^{11}:$

\bigskip

{\bf Lemma 4.1. }{\em  Let  $\delta\in[0,1),$
   $a, \bar b\in S_1^{11},\; W=W(a),\;\bar W=W(\bar b),\; $ and let
$$\mu_1(\delta)=
 {2\over \sqrt
3}\cos\left({\arccos(3\sqrt 3 W)+\pi\over 3 }\right)
-W(1+\delta)\ge$$ $$\mu_2(\delta)=
 {2\over \sqrt
3}\cos\left({\arccos(3\sqrt 3 W)-\pi\over 3 }\right)
-W(1+\delta)\ge$$ $$
  \mu_3(\delta)=
 -{2\over \sqrt
3}\cos\left({\arccos(3\sqrt 3 W) \over 3}\right) -W(1+\delta)
  $$ $($resp.,
  $\bar\mu_1(\delta)\ge\bar\mu_2(\delta)\ge\bar\mu_3(\delta)\;)$
   be the roots of the polynomial
  $$P_{1,\delta}(T,W):=Q_1(T+W+\delta W)=$$
  $$T^3+3W(1+\delta)T^2+(3W^2(1+\delta)^2-1)T+W(1-\delta)+W^3(1+\delta)^3$$
$($resp. of the polynomial
$$ \bar P_{1,\delta}(T,\bar W):=Q_1(T+\bar
W+\delta\bar W)=$$
$$T^3+3\bar W(1+\delta)T^2+(3 \bar W^2(1+\delta)^2-1)T+\bar W(1-\delta)+\bar
 W^3(1+\delta)^3\; ).$$ Then  for any $K>0$ verifying
$ |K-1|+|\bar W-W|\neq 0$ one has
$${1-\delta\over 5+\delta}=:\varepsilon\le { \mu_+(K)\over -\mu_-(K)}\le
 {1\over\varepsilon}= {5+\delta\over 1-\delta}$$
where
$$\mu_-(K):=  \min\{\mu_1(\delta)-K\bar\mu_1(\delta),\;
\mu_2(\delta)-K\bar\mu_2(\delta),\;
\mu_3(\delta)-K\bar\mu_3(\delta)\} ,$$
$$\mu_+(K):=  \max\{\mu_1(\delta)-K\bar\mu_1(\delta),\;
\mu_2(\delta)-K\bar\mu_2(\delta),\;
\mu_3(\delta)-K\bar\mu_3(\delta)\}\;.  $$   }

{\em Proof of Lemma 4.1.} In the proof we will repeatedly use the
following elementary fact:
\medskip

{\em Claim. Let $ l_1\ge l_2 \ge l_3,$ $ l_1+ l_2+ l_3=t\ge 0,$  $
l_3\le -ht,$ with $h>0.$ Then $-l_1/l_3\in [h/(2h+1),(2h+1)/h] $
for $t>0$, $-l_1/l_3\in [1/2,2]$ for $t=0$.}
\medskip

 If $W= \bar W , K=1$ there is nothing to prove. If
$K=1$ one can suppose that $W>\bar W; $ we have
  $$(\mu_1(\delta)-K\bar\mu_1(\delta))+
(\mu_2(\delta)-K\bar\mu_2(\delta))+
(\mu_3(\delta)-K\bar\mu_3(\delta))=3(\bar W -W)(1+\delta)$$ and
$$ \mu_2(\delta)-K\bar\mu_2(\delta)=
 {2\over \sqrt 3}\left(\cos\left({\arccos(3\sqrt 3 W)+ \pi\over 3 }\right )-
 \cos\left({\arccos(3\sqrt 3 \bar W)+\pi\over 3 }\right)\right)$$
$$-(  W -\bar W)(1+\delta)\ge (1-\delta)(  W-\bar W).$$
  Therefore, one can take $\varepsilon=(1-\delta)/(5+\delta) $
  in this case.
 We can suppose then $W>\bar W, K\neq 1. $ Using the
relations $$\mu_1(\delta)(-W)=-\mu_3(\delta)(W),\;
\mu_2(\delta)(-W)=-\mu_2(\delta)(W),\;
\mu_3(\delta)(-W)=-\mu_1(\delta)(W)$$  we can suppose without loss
that $K<1.$

We distinguish then three cases corresponding to different signs
of $W-K \bar W.$  If  $W-K \bar W=0$ then one can take
$\varepsilon=1/2$ since
$$(\mu_1(\delta)-K\bar\mu_1(\delta))+
(\mu_2(\delta)-K\bar\mu_2(\delta))+
(\mu_3(\delta)-K\bar\mu_3(\delta))=0.$$

Let $W-K \bar W=W-\bar W+(1-K)\bar W<0.$ Then
$$(\mu_1(\delta)-K\bar\mu_1(\delta))+
(\mu_2(\delta)-K\bar\mu_2(\delta))+
(\mu_3(\delta)-K\bar\mu_3(\delta))=-3(W-K \bar W)(1+\delta)>0$$
and
$$
   \mu_3(\delta)-K\bar\mu_3(\delta)= \mu_3(\delta)-\bar\mu_3(\delta)
   +(1-K)\bar\mu_3(\delta)=\mu'_3(\delta)(W')(W-\bar W)+(1-K)\bar\mu_3(\delta) $$
for $W'\in (W, \bar W).$ Since
$$\bar\mu_3(\delta)\le {\delta -2\over 3\sqrt 3} < {-1\over 3 \sqrt 3}
\le -\bar W,\;\;\mu'_3(\delta)(W')\le -5/3-\delta\le  -5/3 <-1$$
we get
$$\mu_3(\delta)-K\bar\mu_3(\delta)< -(W-\bar W+(1-K)\bar W)=
-(W-K\bar W)$$ and one can take
$\varepsilon=(2+(3+3\delta))^{-1}=1/(5+3\delta)\ge
(1-\delta)/(5+\delta) $.

Let then $W-K \bar W=W-\bar W+(1-K)\bar W>0.$ We get
$$(\mu_1(\delta)-K\bar\mu_1(\delta))+
(\mu_2(\delta)-K\bar\mu_2(\delta))+
(\mu_3(\delta)-K\bar\mu_3(\delta))=-3(W-K \bar W)<0.$$ If $\bar
W\ge 0 $ then
$$\mu_2(\delta)-K\bar\mu_2(\delta)= \mu_2(\delta)-\bar\mu_2(\delta)
   +(1-K)\bar\mu_2(\delta)=\mu'_2(\delta)(W')(W-\bar W)+(1-K)
   \bar\mu_2(\delta)\ge $$ $$(1-\delta)(W-\bar W) +
   (1-K)(1-\delta)\bar W\ge (1-\delta)(W-K\bar W)$$
   which gives again $\varepsilon=(1-\delta)/(5+\delta).$

    Let $\bar W< 0,\; W\ge 0. $  Then
    $$\mu_2(\delta)-K\bar\mu_2(\delta)\ge (1-\delta) W +K(1-\delta)\bar
    W=(1-\delta)(W-K\bar W).$$
 Let finally $\bar W< 0,\; W< 0. $  Then the same inequality holds
 since the function $f(W):= \mu_2(\delta)(W)/W$ is decreasing for
 $W\in[{-1\over 3\sqrt3},0] $ and $f(0)=(1-\delta). $

\bigskip
This result can be applied to our situation thanks to the
following formulas generalizing those of Section 3; the proofs
remain essentially the same as for Lemma 3.2 (i.e. brute force MAPLE calculation 
together with invariance properties of $w$). Namely, the  matrix $A_{\delta}=H_{\delta}(a)$ 
becomes a block matrix
  $$A_{\delta}=
  \left(%
\begin{array}{cc}
  A_{6,\delta}&0  \\
   0&M_{6,\delta}  \\

\end{array}%
\right)$$ where $A_{6,\delta}=D^2w_{6,\delta}(a')$ is the Hessian of the function 
$$
w_6(a')=P_6(a')/|a'|^{1+\delta}={r_0s_0t_0-r_0s_1t_1-r_1s_0t_1-r_1s_1t_0\over (r_0^2+s_0^2+t_0^2+r_1^2+s_1^2+t_1^2)^{{1+\delta\over 2}}}$$ 
 and $M_{6,\delta}=N_6-(1+\delta)W\cdot I_{6}.$ 

\bigskip
{\bf Lemma 4.2. } {\em  Let $\delta\in[0,1)$ and let $a=(r,s,t)\in
S_1^{11};$   define
$$
W=W(a)=P(a),\; m=m(a) = |s|  ,\; n=n(a)= |t|  .$$

 Then the characteristic polynomial of the
matrix $A_{\delta}= H_{\delta}(a):=D^2w_{\delta}(a)$ is given by

$$ P_{A,\delta}(T)=P_{1,\delta}(T)^2\cdot P_{2,\delta}(T)$$
where

$P_{1,\delta}(T)=P_{1,\delta}(T,W):=Q_1(T+W+\delta W)=$
  $$T^3+3W(1+\delta)T^2+(3W^2(1+\delta)^2-1)T+W(1-\delta)+W^3(1+\delta)^3\; ;$$
$$P_{2,\delta}(T)=P_{2,\delta}(T,W):=T^6+a_{5,\delta}T^5+a_{4,\delta} T^4+a_{3,\delta}T^3+
a_{2,\delta}T^2+a_{1,\delta} T +a_{0,\delta}$$ where
$$a_{5,\delta}:=W(\delta+1)(9-\delta),$$
$$a_{4,\delta}:=W^2(\delta+1)(21+28\delta-5\delta^2)+L(\delta+1)(3-\delta)-2,$$
$$a_{3,\delta}:=-2W(1+\delta)\cdot \left( W^2(\delta+1)(5\delta^2-26\delta-7)-L(2\delta+1)
(3-\delta)+4\right),$$
$$a_{2,\delta}:=-W^4(10\delta^2-53\delta+9)(\delta+1)^3
-2W^2(\delta+1)(3L\delta^3-6L\delta^2-9L\delta +7\delta+3)$$
$$+L\delta^2-3M\delta^2-2L\delta+6M\delta-3L+9M+1$$
$$a_{1,\delta}:=
-(\delta+1) \left(W^4 (5 \delta-3) (\delta-5) (\delta+1)^3-2
(\delta+1) (-2 L \delta^3 +5 L \delta^2 +4 L \delta-6 \delta -3
L+2) W^2\right.$$
$$\left.+2 (3-\delta) (-3 \delta M+L \delta-L)+1-\delta\right)W$$
$$ a_{0,\delta} =(1-\delta) \left(W^6 (\delta-5) (\delta+1)^5+
W^4 (\delta+1)^3 (L \delta^2-2 L \delta-3 L+4)\right.$$
$$\left.-W^2 (\delta+1) (L \delta^2-3 M \delta^2+\delta+6 M
\delta-4 L \delta-1+3 L+9 M)-M (1-\delta)\right)\; , $$ with
$L=m^2+n^2-n^2m^2-n^4-m^4,\;  M=m^2n^2(1-n^2-m^2)$ as before}.

\medskip
 A MAPLE calculation gives then  for the resultant

$$R_{\delta}(r_0,m,n):=Res(P_{1,\delta},P_{2,\delta})=
16m^4n^4(1-n^2-m^2-r_0^2)^3\cdot R(W,\delta)$$ where
 $$R(W,\delta)=27(\delta+1)^3(3-\delta)^3W^4+9(\delta-1)^2(\delta-3)^2(\delta+1)^2
W^2-(\delta-1)^2(\delta^2-2\delta-2)^2.$$  Denote by
$W_0(\delta)\in (0,1/3\sqrt 3]$ the unique positive root of
$R(W,\delta).$  Recall that the set $\Phi(S^{11}) $ of possible
triples $\Phi(a)=(r_0,m,n): r_0=r_0(a), m=m(a),n=n(a)$ for $a\in
S^{11}_1 $ is a quarter $\bar B_{++}:=B_1\bigcap \{ m\ge 0, n\ge
0\} $ of the closed unit ball $B=B_1 \subset V;$ recall also that
$ W(a)=r_0mn$. Let $ B_+(\delta)$ (resp. $ B_-(\delta)$, $
B_0(\delta)$) be the subset of $(r_0,m,n)\in \Phi(S^{11})$ where
$R_{\delta}(W)>0$ (resp. $R_{\delta}(W)<0$, $R_{\delta}(W)=0$).
Then

\medskip
 $ B_0(\delta)=S_{++}^2\bigcup D_{r_{0}+}\bigcup D_{m+}\bigcup
D_{n+}$ with $D_{m+}=\bar B_{++}\bigcap \{ m=0\}$ etc.,
$$ B_+(\delta)=B_{++}\bigcap \{ r_0mn  >W_0(\delta)\},\;
\bar B_+(\delta)=\bar B_{++}\bigcap \{ r_0mn  \ge W_0(\delta)\},$$
$$ B_-(\delta)=B_{++}\bigcap \{0< r_0mn  <W_0(\delta)\},\;
\bar B_-(\delta)=\bar B_{++}\bigcap \{ r_0mn  \le W_0(\delta)\}.$$

Note that these sets are invariant under the reflection $Refl:
(r_0,m,n)\mapsto (-r_0,m,n) $; $ B^0(\delta)$ and $\bar
B_-(\delta)$ are connected, while $ B_-(\delta)$, $\bar
B_+(\delta)$ and $B_+(\delta)$ have two connected components each.

\bigskip {\bf Lemma 4.3. }
 {\em  Let $a \in S_1^{11},$ let
 $\lambda_{1}(\delta,a)\ge\lambda_{2}(\delta,a)\ge\ldots\ge\lambda_{12}(\delta,a)$
 be the eigenvalues of $D^2w_{\delta}(a)$ and let
 $\mu_1(\delta,a)\ge\mu_2(\delta,a)\ge\mu_3(\delta,a)$ be the roots of
 $P_{1,\delta}(T,W(a)).$ Then

 $$(i)\;\;\lambda_{1}(\delta,a)=\lambda_{2}(\delta,a)=\mu_1(\delta,a),\;
 \lambda_{12}(\delta,a)=\lambda_{11}(\delta,a)=\mu_3(\delta,a)\; ; \hskip 1. cm$$
$$(ii)\;\;\lambda_{5}(\delta,a)=\lambda_{6}(\delta,a)=\mu_2(\delta,a)\;\;
 {\hbox {for}}\; \;
 \Phi(a)\in \bar B_+(\delta),\;W=W(a)\ge 0\; ;\quad$$
 $$(iii)\;\;\lambda_{7}(\delta,a)=\lambda_{8}(\delta,a)=\mu_2(\delta,a)\;\;
 {\hbox {for}}\; \;
 \Phi(a)\in \bar B_+(\delta),\;W=W(a)\le 0\;;\quad$$
 $$(iv)\;\;\lambda_{6}(\delta,a)=\lambda_{7}(\delta,a)=\mu_2(\delta,a)\;\;
 {\hbox {for}}\; \;
 \Phi(a)\in \bar B_-(\delta)\; .\hskip 2.9 cm$$}

 \medskip
{\em Proof of Lemma 4.3.} Since
$\lambda_{1}(\delta,a)=\lambda_{12}(\delta,-a),$
$\lambda_{12}(\delta,a)=\lambda_{1}(\delta,-a),$
$\lambda_{6}(\delta,a)=\lambda_{7}(\delta,-a),$
$\lambda_{7}(\delta,a)=\lambda_{6}(\delta,-a),$
$\lambda_{8}(\delta,a)=\lambda_{5}(\delta,-a),$
$\lambda_{5}(\delta,a)=\lambda_{8}(\delta,-a),$
 $ W(-a)=-W(a),$ $
(iii)$ is implied by  $(ii)$ and, moreover one can suppose without
loss that $\Phi(a)=(r_0,m,n)\in {\bf R}_+^3.$ Since in the
interior of the domain $B_+(\delta)\bigcap{\bf R}_+^3$ (resp.
$B_-(\delta) \bigcap{\bf R}_+^3$) the function
$R_{\delta}(r_0,m,n)$ does not vanish, it is sufficient to verify
the ordering of the roots at a single point in
$B_-(\delta)\bigcap{\bf R}_+^3$ (resp. at a single point in
$B_+(\delta) \bigcap{\bf R}_+^3).$  We use
$a_-:=(\varepsilon,\varepsilon,\varepsilon)\in
B_-(\delta)\bigcap{\bf R}_+^3$ and  $a_+:=(1/\sqrt 3,1/\sqrt
3,1/\sqrt 3-\varepsilon)\in B_+(\delta)\bigcap{\bf R}_+^3$ for
sufficiently small $\varepsilon>0.$  Let
$\nu_1(\delta,a)\ge\nu_2(\delta,a)\ge\ldots \ge\nu_6(\delta,a) $
be the roots of $P_{2,\delta}(T,W(a)).$ Elementary calculations
show that for $a=a_-$ one has $W=W(a)=\varepsilon^3,$
$$\mu_1(\delta,a)=1+O(\varepsilon^3),\;
\mu_2(\delta,a)=O(\varepsilon^3),\;\mu_3(\delta,a)=-1+O(\varepsilon^3),$$
while $P_{2,\delta}(T,W(a))=F_1(T,\varepsilon)\cdot
F_2(T,\varepsilon)$

where $$F_1(T,\varepsilon)=   T^2-1+2
\varepsilon^2+O(\varepsilon^3), $$
$$F_2(T,\varepsilon)=T^4+\varepsilon^3 (7 +6 \delta -\delta^2
)T^3$$
 $$ +(12 \varepsilon^6 \delta^2 +3 \varepsilon^4 \delta^2 -3
\varepsilon^6 \delta^3 +21 \varepsilon^6 \delta+6 \varepsilon^6 +4
\varepsilon^2 -2 \varepsilon^2 \delta^2 -1-6 \varepsilon^4 \delta
+4 \varepsilon^2 \delta -9 \varepsilon^4 )T^2$$
 $$ +\varepsilon^3(1-10
\varepsilon^2 - \delta^2 -12 \varepsilon^4 \delta^2   +4
\varepsilon^2 \delta -18 \varepsilon^4 \delta +10 \varepsilon^2
\delta^2 -4 \varepsilon^2 \delta^3 +6 \varepsilon^4 \delta^3+
O(\varepsilon^6))T $$
 $$ +
\varepsilon^4(1-\delta)^2-\varepsilon^6(\delta+1)(\delta-1)^2
(2\delta\varepsilon^2-4\varepsilon^2+1) +
 O(\varepsilon^{10})$$
 and thus
 $$\mu_1(\delta,a)\ge\nu_1(\delta,a)=1-\varepsilon^2+O(\varepsilon^3)\ge
\nu_2(\delta,a)=1-\varepsilon^2(2+2\delta-\delta^2)+O(\varepsilon^3),$$
$$ \nu_3(\delta,a)=(1-\delta)\varepsilon^2 +O(\varepsilon^3)\ge\mu_2(\delta,a)\ge
\nu_4(\delta,a)=-(1-\delta)\varepsilon^2+O(\varepsilon^3),$$
$$
\nu_5(\delta,a)=-1+\varepsilon^2(2+2\delta-\delta^2)+O(\varepsilon^3)\ge
\nu_6(\delta,a)=-1+\varepsilon^2+O(\varepsilon^3)\ge\mu_3(\delta,a)
$$
which proves the claim in this case.

 \medskip
For $a=a_+$ one has $W=W(a)=(1/\sqrt 3-\varepsilon)/3$ and similar
calculations give
$$\mu_1(\delta,a)={2-\delta\over 3\sqrt 3 }+3^{-1/4} \sqrt{2\varepsilon}
+O(\varepsilon),\; \mu_2(\delta,a)={2-\delta\over 3\sqrt 3 }-
3^{-1/4} \sqrt{2\varepsilon} +O(\varepsilon),\;$$ $$
\mu_3(\delta,a)={-7-\delta\over 3\sqrt 3
}+(5/3+\delta)\varepsilon+O(\varepsilon^2),$$

while $P_{2,\delta}(T,W(a))=G_1(T,\varepsilon)\cdot
G_2(T,\varepsilon)^2\cdot G_3(T,\varepsilon)^2$

where
$$G_1(T,\varepsilon):=T^2+{(\delta+1)\over 3\sqrt 3}(5-\delta)
(1-\sqrt 3\varepsilon)T $$
$$+{(1-\delta)\over 27}(3\delta^2\varepsilon^2-2\sqrt
3\delta^2\varepsilon-12\delta\varepsilon^2+\delta^2 +8\sqrt
3\delta\varepsilon-15\varepsilon^2+5\delta+10\sqrt 3
\varepsilon-14)\;,$$

$$G_2(T,\varepsilon):= T+{4+\delta\over 3\sqrt 3} -{\varepsilon(\delta+1)\over 3}
\; ,$$
$$G_3(T,\varepsilon):= T-{2-\delta\over 3\sqrt 3} -{\varepsilon(\delta+1)\over 3}
  \; ,$$
and thus
$$\mu_1(\delta,a)\ge\nu_1(\delta,a)=
\nu_2(\delta,a)= {2-\delta\over 3\sqrt 3} +O(\varepsilon)\ge
\mu_2(\delta,a)\ge$$
$$  \nu_3(\delta,a)= {(2-\delta)(1-\delta)\over 3\sqrt 3}
+O(\varepsilon)\ge \nu_4(\delta,a)=\nu_5(\delta,a)=-{4+\delta\over
3\sqrt 3} +O(\varepsilon)\; ,$$
$$
\nu_6(\delta,a)={-7-\delta\over 3\sqrt 3
}+{\varepsilon(\delta-5)(\delta-9)(\delta+1)\over
3(9-2\delta+\delta^2)}+O(\varepsilon^2) \ge\mu_3(\delta,a)
$$
which finishes the proof of the lemma (note that
$$\nu_6(\delta,a)-\mu_3(\delta,a)=
{2\varepsilon\delta(7+\delta)(1-\delta)\over
3(9-2\delta+\delta^2)}+O(\varepsilon^2)\ge 0\; ).$$

\bigskip
{\em End of  proof of Proposition 4.1.} If $W(a)$ and $W(b)$ are
of the same sign we get the result applying Lemmas 4.1 and 4.3
with $K:= |b| ^{-\delta} $; in the exceptional case $K=1, \;
W(a)=W(b)$ the trace of $H_{\delta}(a,b,O)$ vanishes and the claim
is valid for $C_{\delta}=11$. In the case $W(a)\cdot W(b)<0$ we
can suppose without loss that $W(a)>0,\; W(b)<0;$ if $\Phi(a)
\notin B_+$ or $\Phi(\bar b) \notin B_+$ then Lemmas 4.1 and 4.3
work as well. Thus we can suppose $\Phi(a) \in B_+,\;\Phi(\bar b)
\in B_+;$ then
$$Refl(\Phi(\bar b)) \in B_+,\;W(-\bar b)>0,\;
\lambda_{i}(-b)=-\lambda_{13-i}(b),\;\lambda_{i}(-\bar
b)=-\lambda_{13-i}(\bar b)\;$$ and
$$
Tr(H_{\delta}(a,b,O))=-(W(a)+KW(-\bar b))(\delta+1)(15-\delta)<0$$
which implies immediately that $11\ge -\Lambda_{1}/\Lambda_{12}. $

Moreover, $$\Lambda_{1}\ge  \lambda_{6}(a)-K \lambda_{6}(\bar
b)=\lambda_{6}(a)+K \lambda_{7}(-\bar b)=
\mu_2(\delta,a)+K\mu_2(\delta, -\bar b)\ge $$ $$
(1-\delta)(W(a)+KW(-\bar b))= {(1-\delta)Tr(H_{\delta}(a,b,O))
\over (\delta+1)(15-\delta)}>0$$ and thus $$
-\Lambda_{1}/\Lambda_{12}\ge \left(11+ {
(\delta+1)(15-\delta)\over 1-\delta}\right)^{-1}={ 1-\delta \over
26+3\delta-\delta^2}$$ which finishes the proof of the
proposition.

\bigskip To deduce Corollary 1.1 we need the map
$$H_{\delta}: B^{12}_1-\{0\} \longrightarrow Q\;,\;\; a \mapsto D^2w_{\delta}(a)$$
where $Q=S^2({\bf R}^{12})$ denotes the space of quadratic forms
on ${\bf R}^{12}.$ The following result is sufficient  to conclude
using Proposition 4.1 and Lemma 2.2 of [NV1]:

\bigskip
{\bf Lemma 4.4.} {\em Let $\delta \in (0,1).$ Then the image
$H_{\delta}\left (B^{12}_1-\{0\}\right )\subset Q$ is
diffeomorphic to the product $V_{11,\delta}\times [1,\infty)$ with
a smooth 11-dimensional manifold $V_{11,\delta}$.}

\bigskip
{\em Proof.} Since
$D^2w_{\delta}(a)=D^2w_{\delta}(a/|a|)|a|^{-\delta}$ it is
sufficient to show two facts:

\medskip (i)  $H_{{\delta}\;| S^{11}_1}:S^{11}_1 \longrightarrow Q$
is a smooth embedding;

\medskip
(ii) if $D^2w_{\delta}(a)=D^2w_{\delta}(b)\cdot k$ with $k>0$ then
$k=1$.

\medskip
Lemmas 4.1 and 4.2 imply (ii). To prove (i) we fix
 $a\neq b \in S^{11}_1 $ and
consider $ d= {a-b\over  |a-b|}\in S^{11}_1.$ Let then $e, f\in
S_1^{11}\bigcap a^{\perp} \bigcap  b^{\perp}. $ Since $e, f \perp
a, b $ one has
$$w_{\delta,ee}(a)=P_{ee}(a)-(1+\delta)P(a),\;\;
w_{\delta,ee}(b)=P_{ee}(b)-(1+\delta)P(b),$$
$$w_{\delta,ff}(a)=P_{ff}(a)-(1+\delta)P(a),\;\;
w_{\delta,ff}(b)=P_{ff}(b)-(1+\delta)P(b)$$ and hence
$$\left(w_{\delta,ee}(a)-w_{\delta,ee}(b)\right)-
\left(w_{\delta,ff}(a)-w_{\delta,ff}(b)\right)=$$
$$\left(P_{ee}(a)-P_{ee}(b)\right)
-\left(P_{ff}(a)-P_{ff}(b)\right)= |a-b| (P_{eed}-P_{ffd})\ge
{2\over \sqrt 3} |a-b| $$ for suitable vectors $e,f$ as in the
proof of Proposition 2 in [NV1, Section 4]. It follows that $$
\max\{|w_{\delta,ee}(a)-w_{\delta,ee}(b)|,|w_{\delta,ff}(a)-w_{\delta,ff}(b)
|\}\ge  |a-b|/\sqrt 3 $$ which finishes the proof.
\section{Isaacs equation}

\medskip
We can then prove Theorem 1.2 as a simple consequence of the results of Section 4.
Denote by $K_C\subset S^2({\bf R}^2) $ the cone of positive symmetric matrix with
the ellipticity constant $C$, i.e., if $A\in K_C, \  A=\{ a_{ij} \}$ then
$$C^{-1}|\xi  |^2 \leq \sum a_{ij}\xi_i\xi_j \leq C|\xi |^2.$$
Recall the following results from [NV3, Section 5]:

\medskip
{\bf Lemma 5.1.} {\it Let $w\in C^{\infty }({\bf R}^n\setminus 0) $ be a homogeneous order $\alpha , 1<\alpha
\leq 2$ function. Assume that for any two points $x,y\in {\bf R}^n,\ 0<|x|,|y|\le 1$, there exists a matrix
$A\in K_C$ orthogonal to both forms $D^2w(x), D^2w(y),$ 
$$Tr(A D^2w(x))=Tr (AD^2w(y))=0.$$
Then $w$ is a viscosity solution to an Isaacs equation. }

\medskip
 Recall that a symmetric matrix $A$ is  called strictly hyperbolic if
$$\frac{1}{ M} < -\frac{\lambda_1(A)}{\lambda_n (A)} < M$$
for a positive $M$.

\medskip
{\bf Lemma 5.2.} {\it Let $F_1, F_2$ be two quadratic forms  in ${\bf R}^n$  s.t. the  form  $ \alpha F_1+\beta  F_2$
 is strictly hyperbolic for any $ (\alpha, \beta) \in {\bf R}^2\setminus \{0\}$. Then there exists a positive   quadratic form
 $Q$ orthogonal to both forms $F_1, F_2$, 
$$Tr(F_1Q)=Tr(F_2Q)=0. $$ }
\medskip
 The results of Section 4 imply that the 
form  $ \alpha {D^2w_{\delta}} _{|H}(x) -\beta {D^2w_{\delta}} _{|H}(y) $
 is strictly hyperbolic for positive $ \alpha, \beta$; since the function $w_{\delta}$ is 
odd, it remains true for any  $ (\alpha, \beta) \in {\bf R}^2\setminus \{0\}$  and thus Lemmas 5.1 and
5.2 imply Theorem 1.4.

\bigskip

\section{ Eleven Dimensions}

For a unit vector  $a\in S_1^{10} \subset  {\bf R}^{11}$ we
continue to denote $D^2 w_H(a)$ by $H(a)$.

\bigskip
{\bf Lemma 6.1. }{\em Let $a \in S_1^{10} $ and let $\lambda_1\ge
\lambda_2\ge\ldots\ge\lambda_{11} $ be the eigenvalues  of
$A=H(a)$. Then}
$$\lambda_6 ={2\over \sqrt 3}\cos\left({\arccos(3\sqrt 3 P_{H}(a))+\pi\over 3
}\right)-P_{H}(a).$$

\medskip
{\em Proof.} This follows from Lemma 3.1 and Lemma 2.3.

\medskip
 Let then   $a\neq b\in S_1^{10} $  and let $O\in {\hbox {O}}({11} )$
be an orthogonal matrix s.t. $H(a,b,O):=H(a)- {^tO}\cdot H(b)\cdot
O\neq 0$. Denote $ \Lambda_1\ge\Lambda_2\ge \ldots\ge\Lambda_{11}$
the eigenvalues of the matrix $H(a,b,O).$ As in Section 3 above one gets

\bigskip
{\bf Lemma 6.2. } {\em  Let $ A:=H(a),$ $B:={^tO}\cdot H(b)\cdot
O.$

$(i)$ If $P_{H}(a)-P_{H}(b)\ge 0 $ then }${\hbox
{Tr}}(B-A)=14(P_{H}(a)-P_{H}(b))\le 14\Lambda_1; $

$(ii)$ {\em If $P_{H}(a)-P_{H}(b)\le 0 $ then} ${\hbox
{Tr}}(B-A)=14(P_{H}(a)-P_{H}(b))\ge 14\Lambda_{11}. $

 which implies

\bigskip
 {\bf Proposition 6.1.} {\em Let  $a\neq b\in S_1^{10} $  and
 let $O\in {\hbox {O}}({11} )$ be an orthogonal matrix s.t.
 $H(a,b,O):=H(a)- {^tO}\cdot H(b)\cdot O\neq 0$.
 Denote $ \Lambda_1\ge\Lambda_2\ge
 \ldots\ge\Lambda_{11}$  the eigenvalues of the matrix
 $H(a,b,O).$
  Then}

$${1\over 24} \le {\Lambda_1\over -\Lambda_{12}}\le 24.$$
\medskip
 
 Proposition 6.1 and Lemma 2.1 give a proof of Theorem 1.3 exactly  
as Proposition 3.1 implies Theorem 1.1 in the case $\delta=0$.  

\medskip
{\em Remark 6.1.}  The ellipticity constant $C$ of thus obtained
functional $F$ verifies
 $C\le 4\cdot 24^2\sqrt{11}<10^4. $  

\medskip
 {\em Remark 6.2.} One can not directly use the approach of   Section 4 to the 
function $$w_H/ |x|^{\delta } $$ for $\delta >0$ since
 although the corresponding Hessian $D^2(w/ |x|^{\delta }) $ always
 has double  eigenvalues,  they  position in the spectrum is
 not fixed and can vary from (5,6) to (7,8), see Lemma 4.3 above.
 It means  that after the restriction on a hyperplane $H$ we
 lose the property necessary to control the ellipticity and thus can not construct a singular solution 
in 11  dimensions.

  \section{Singular solutions with cusp}

  Let $P$ be a linear elliptic operator of the form
$$P = \sum_{i,j} a_{ij}(x) {\partial^2 \over \partial x_i\partial x_j },$$
defined in a half-ball $B_+=\{ x\in B\subset {\bf R}^n, x_1>0\} $,
 $a_{ij} \in L_{ \infty} (B_+)$ and satisfying the inequalities
 $$C^{-1}|\xi|^2\le \sum
a_{ij}\xi_i\xi_j\le C |\xi |^2\;,\forall\xi\in {\bf R}^n\;.$$

Let $z\in C^2(B_+)$  and $Pz=0$ in   $B_+, $ $z=0$  on $L,$ where
$L=\{ x\in B, x_1=0 \}$. Assume that $z<1$ in $B_+$. Then it is
well known, [GT], that
$$ |\nabla z(0)| \leq K,$$
where constant $K$ depends on the ellipticity constant $C$.

\bigskip
{\bf Lemma 7.1.} {\em The following inequality holds with
positive constants $K, \epsilon $ depending on the ellipticity
constant $C$:
$$ |z-dz(0)|\leq K|x|^{1+\epsilon },$$
where $dz$ is the differential of the function $z$}.

\medskip The lemma follows directly from P. Bauman's boundary
Harnack inequality, [B].

\bigskip
{\em Proof of Theorem 1.4. } We may assume w.r.g.  that $F(0)=0$,
otherwise instead of the function $u$ we consider the function
$u+c|x|^2$ with a suitable constant $c$.

Set
$$v(r)=\sup_{|x|=r} u(x),$$
$$u_i=u(x_1,...,-x_i,...,x_n),$$
$$z_i=u-u_i.$$
Since $u$ is a solution of a Hessian equation the functions $u_i$
are solutions of the same equation as well. Hence functions $z_i$
given as the difference of two solutions of the fully nonlinear
elliptic equation are solutions to a linear elliptic equation
$Pz_i=0$ in $B$.  Define a linear function $l$ as
$$l={1\over 2} \sum dz_i(0).$$

Set
$$u_0=u-l.$$

Let $|y|=|y'|=r<1$. Choose in ${\bf R}^n$ an orthonormal
coordinate system $y_1,...,y_n$, such that $y_1=(y-y')/|y-y'|$.
Set
$$u'(y_1,...,y_n)=u_0(-y_1,...,y_n),$$
 $$v=u_0-u'.$$

 Since $F(u')=0$ we get $Pv=0$ in $B$. Moreover
 $$\nabla v(0)=0. $$

 Hence by  Lemma 7.1,
 $$v(x)=o(|x|^{1+\epsilon }).$$

 Therefore
 $$u_0(y)-u_0(y')= o(|y|^{1+\epsilon }).$$

 Set
 $$h(r) = \inf_{|x|=r} u_0(x),$$
 $$h_0(r) = \sup_{|x|=r} u_0(x).$$
 Then
  $$h(|x|)-h_0(|x|)= o(|x|^{1+\epsilon }).\leqno (7.1)$$

  Since $F(0)=0$, we may assume without loss that $u(0)=0, \  h_0'(1)>0$.
Then  by the maximum principle $h_0'(r)$ is a monotone function of
$r$. If $h(r)= o(|x|^{1+\epsilon /2})$ we may set $h\equiv 0$
  and the theorem is proved. Assume that $h(r) >\epsilon |x|^{1+\epsilon /2}$. Then from
  (6.1) it follows that $|h(r)|$ is a positive function for sufficiently small $r$.

  By a direct computation
  $$\lambda (D^2 h(|x|))= (h'' , h'/|x|, ...,h'/|x|).$$

  Hence $h$ has no local minimums and since $h>0$ we get $h'>0, h''<0$ for sufficiently small $r$.
  Therefore $h$ is a monotone, concave function for small $r$.

  For any $0<r<1$ there exists a point $x_0, |x_0|=r$ such that $u_0(x_0)=h(r)$ and since
$h-u_0 \leq 0$ the quadratic part of the function $u_0-h$ is
non-negatively defined. Hence from the uniform ellipticity
condition for $F$ we get the inequality
$$-|x|h''/h'>\delta ,$$
on an interval $(0,a)$ for some $a>0$, where $\delta $ depends on
the ellipticity constant. From the last inequality it follows that
$$h(r)>r^{1-\delta }$$
on $(0,a)$. Since we can redefine $h$ on $(a,1)$ as a monotone,
concave function, the theorem is proved.

\pagebreak
 \centerline{\bf REFERENCES}

\medskip\medskip
\medskip

\noindent [A] A.D. Alexandroff, {\it Sur les th\'eor\`emes
d'unicite pour les surfaces ferm\'ees}, Dokl. Acad. Nauk 22
(1939), 99--102.

\medskip
\noindent [B] P. Bauman, {\it Positive solutions of elliptic
equations in non-divergence form and their adjoints}, Ark. Mat. 22
(1984), 153--173.

\medskip
\noindent [Ba] J. Ball, {\it Convexity conditions and existence
theorems in nonlinear elasticity}, Arch. Rat. Mech. Anal. 63
(1977), 337--403.

\medskip
 \noindent [C] L. Caffarelli,  {\it Interior a priory estimates for solutions
 of fully nonlinear equations}, Ann. Math. 130 (1989), 189--213.

\medskip
 \noindent [CC] L. Caffarelli, X. Cabre, {\it Fully Nonlinear Elliptic
Equations}, Amer. Math. Soc., Providence, R.I., 1995.

\medskip
 \noindent [CIL]  M.G. Crandall, H. Ishii, P-L. Lions, {\it User's
guide to viscosity solutions of second order partial differential
equations,} Bull. Amer. Math. Soc. (N.S.), 27(1) (1992), 1--67.

\medskip
 \noindent [CNS] L. Caffarelli, L. Nirenberg, J. Spruck, {\it The Dirichlet
 problem for nonlinear second order elliptic equations III. Functions
  of the eigenvalues of the Hessian, } Acta Math.
   155 (1985), no. 3-4, 261--301.

\medskip
 \noindent [E] L. C. Evans, {\it Classical solutions of fully nonlinear,
 convex, second-order elliptic equations }, Comm. Pure Appl. Math.
 35 (1982), 333--363.

\medskip
 \noindent [Fu]  W. Fulton, {\it Eigenvalues, invariant factors, highest weights,
and Schubert calculus,} Bull. Amer. Math. Soc. (N.S.), 37(3)
(2000), 209--249.

\medskip
 \noindent [GT] D. Gilbarg, N. Trudinger, {\it Elliptic Partial
Differential Equations of Second Order, 2nd ed.}, Springer-Verlag,
Berlin-Heidelberg-New York-Tokyo, 1983.

\medskip
 \noindent [HNY] Q. Han, N. Nadirashvili, Y. Yuan, {\it Linearity of
homogeneous order-one solutions to elliptic equations in dimension
three,} Comm. Pure Appl. Math. 56 (2003), 425--432.

\medskip
 \noindent [HL]
R. Harvey, H. B. Lawson Jr., {\it Calibrated geometries,}
Acta Math. 148 (1982), 47--157.

\medskip
 \noindent [JX] J. Jost., Y.-L. Xin, {\it A Bernstein theorem
for special Lagrangian graphs,} Calc. Var. Part. Diff. Eq. 15
(2002), 299--312.

\medskip

\noindent [K] N.V. Krylov, {\it Nonlinear Elliptic and Parabolic
Equations of Second Order}, Reidel, 1987.

\medskip
\noindent [LO] H. B. Lawson Jr., R. Osserman,
{\it Non-existence, non-uniqueness and irregularity of solutions to the minimal surface system,}
Acta Math. 139 (1977),  1--17.

\medskip
\noindent [N] L. Nirenberg, {\it On nonlinear elliptic partial
differential equations and H\"older continuity,} Comm. Pure Appl.
Math. 6 (1953), 103--156.

\medskip
\noindent [NV1] N. Nadirashvili, S. Vl\u adu\c t, {\it
Nonclassical solutions of fully nonlinear elliptic equations,}
Geom. Func. An. 17 (2007), 1283--1296.

\medskip
\noindent [NV2] N. Nadirashvili, S. Vl\u adu\c t, {\it Singular
solutions to fully nonlinear elliptic equations,} J. Math. Pures
Appl. 89 (2008), 107--113.

\medskip
\noindent [NV3] N. Nadirashvili, S. Vl\u adu\c t, {\it
Nonclassical solutions of fully nonlinear elliptic equations II.
 Hessian  Equations  and Octonions,}  arXiv:0912.3126, submitted.

\medskip
 \noindent [NY] N. Nadirashvili, Y. Yuan, {\it Homogeneous solutions
to fully nonlinear elliptic equation}, Proc. AMS, 134:6 (2006),
1647--1649.

\medskip
\noindent [T1] N. Trudinger, {\it Weak solutions of Hessian
equations,} Comm. Partial Differential Equations 22 (1997), no.
7-8, 1251--1261.

\medskip
\noindent [T2] N. Trudinger, {\it On the Dirichlet problem for
Hessian equations,} Acta Math.
   175 (1995), no. 2, 151--164.

\medskip
\noindent [T3] N. Trudinger, {\it H\"older gradient estimates for fully nonlinear elliptic equations,}
Proc. Roy. Soc. Edinburgh Sect. A 108 (1988), 57--65.

\medskip
\noindent [We] G. Weyl, {\it Das asymptotische Verteilungsgezets
des Eigenwerte lineare partieller Differentialgleichungen,} Math.
Ann.
   71 (1912), no. 2,  441--479.
\end{document}